\documentclass[12pt,reqno]{amsart}
\usepackage[english]{babel}
\usepackage[utf8]{inputenc}
\usepackage[T1]{fontenc}
\usepackage{wrapfig}

\usepackage[top=25mm, bottom=20mm, left=35mm, right=35mm]{geometry}

\usepackage{bbm}
\usepackage{float, verbatim}
\usepackage{amsmath}
\usepackage{amssymb}
\usepackage{amsthm}
\usepackage{amsfonts}
\usepackage{graphicx}
\usepackage{mathtools}
\usepackage[pagebackref]{hyperref}
\usepackage{color}
\usepackage{tikz}
\usepackage{xcolor}
\usepackage{pgfplots}
\usepackage{caption}
\usepackage{subcaption}
\usepackage{enumerate}
\usepackage{url}
\usepackage{needspace}

\hypersetup{
	colorlinks=true,
	linkcolor=blue,
	filecolor=magenta,      
	urlcolor=cyan,
	citecolor= red,
}

\pgfplotsset{compat=1.18}

\newtheorem*{thm*}{Theorem}
\newtheorem*{conj*}{Conjecture}
\newtheorem*{ques*}{Question}
\newtheorem*{rem*}{Remark}
\newtheorem*{defn*}{Definition}
\newtheorem*{mainques*}{Main questions}

\newtheorem{thm}{Theorem}[section]

\newtheorem{lem}[thm]{Lemma}
\newtheorem{cor}[thm]{Corollary}
\newtheorem{defn}[thm]{Definition}

\newtheorem{prop}[thm]{Proposition}
\newtheorem{conj}[thm]{Conjecture}
\newtheorem{rem}[thm]{Remark}

\numberwithin{equation}{section}

\def \one{\mathbf{1}}

\def \bN {{\mathbb N}}

\def \bZ {{\mathbb Z}}

\def \epsilon {{\varepsilon}}

\def \bC {\mathbb C}

\def \bN {\mathbb N}

\def \bZ {\mathbb Z}

\def \cF {\mathcal F}

\def \leq {\leqslant}

\def \geq {\geqslant}

\def \ds1 {\mathds{1}}

\title[Ternary expansions of $2^n$]{Ternary expansions of $2^n$}

\author{Han Yu}
\address{Han Yu, College of Mathematics and Statistics, Center of Mathematics,
Chongqing University, Chongqing, 401331, China}
\email{han.yu.2@cqu.edu.cn}
\date{}

\begin{document}

\begin{abstract}
Consider the set $T$ of integers $n>0$ so that $2^n$ has no digit $2$ in
its ternary expansion. A conjecture of Erd\H{o}s asserts that $T$ is finite.
Narkiewicz's classical result gives $\#(T\cap[1,N])\ll N^{\log 2/\log 3}$.
We show that for some $\delta>0$,
\[
 \#(T\cap[1,N])\ll N^{\log 2/\log 3-\delta}.
\]
This answers a question of Lagarias. The main step is a modification of a Fourier decay argument in Banaji-Yu \cite{BY} for the discrete setting.
\end{abstract}
\maketitle

\section{Introduction}

Let $A\subset\bN$ be an integer sequence. Let $C=C_{3,\{0,1\}}$ be the
set of positive integers whose ternary expansions have only digits $0,1$.
We consider the intersection $A\cap C$. If there are no obvious ternary
digit constraints on the generation of $A$, we should expect that $A$
and $C$ are in some sense independent. For some digit restrictions,
such independence can be established for primes. An example is
Maynard's theorem on primes omitting a prescribed decimal digit
\cite{May2019}.

Write $\alpha=\log 2/\log 3$. The set $C$ has about $X^\alpha$ elements
up to $X$. If one follows the independence heuristic, a sufficiently
thin sequence $A$ should have only finitely many elements in $C$.
For example, if $\#(A\cap[1,X])\ll X^\epsilon$ for some
$0<\epsilon<1-\alpha$, the expected number of elements of $A\cap C$
between $3^k$ and $3^{k+1}$ is at most a constant times
$3^{k(\epsilon+\alpha-1)}$, which is summable in $k$.
The following conjecture is due to Erd\H{o}s \cite{Erd1979}. See also Problem$\#406$ on $www.erdosproblems.com$

\begin{conj}
There are only finitely many integers $n\geq1$ so that $2^n\in C$.
\end{conj}

More precisely, Erd\H{o}s conjectured that the only such positive
exponents are $2$ and $8$. Clearly, the finiteness conjecture fails if
one takes $3^n$ instead of $2^n$. This reveals an obstruction to the
independence heuristic. 

Let
\[
 T=\{n\in\bN:2^n\in C\},\qquad \bN=\{1,2,\dots\}.
\]
The following classical result is due to Narkiewicz \cite{Nar1980};
see also the introduction of Lagarias \cite{Lag2009}.

\begin{thm}[Narkiewicz]
For $N\in\bN$,
\[
 \#(T\cap[1,N])\ll N^\alpha.
\]
\end{thm}
Here $\alpha=\log 2/\log 3$ can be seen as the ``dimension'' of $C$ which naturally serves as the above upper bound which will be illustrate in Section \ref{sec: idea}. We show the following improvement and answers a question of Lagarias \cite[Section 1.6]{Lag2009}.

\begin{thm}\label{thm:main}
Let
\begin{align*}
 \delta&      =\frac{(2\alpha-1)(1-\alpha)}{5-2\alpha}
       =0.0258536438\ldots.
\end{align*}
For every interval $I$ of $N\geq1$ consecutive positive integers,
\[
 \#(T\cap I)\ll N^{\alpha-\delta},
\]
with an absolute implied constant. In particular,
$\#(T\cap[1,N])\ll N^{\alpha-\delta}$.
\end{thm}
Since $1/2<\alpha<1$, we have $\delta>0$. In this paper, we only
present the simple $L^2$ argument leading to this value of $\delta$.
This was our original argument, whose main idea comes from the work of
Banaji--Yu \cite{BY} in the continuous setting.

We are also curious whether a larger $\delta$ can be obtained using
$L^1$ arguments or refined Fourier decay analysis (circle method). For this purpose,
AI was used and the answer is affirmative. But the resulting arguments
are significantly more involved.
We decided to keep the simple argument and not to pursue the
$\delta$ race here.\footnote{Limited by the author's accessible computation resource, there is no sign of $\alpha-\delta$ approaching $0$, which would be the next milestone. Nonetheless, we will make the AI-generated proofs for larger $\delta$ public,
without seeking journal publication.}

\begin{rem}
    The key method for proving Theorem \ref{thm:main} is an $L^2$-argument for Fourier decay (Proposition \ref{prop:l2}) for non-linear images of self-similar measures in the discrete setting. It is of course possible to apply this method for: 1. other discrete self-similar sets and 2. other non-linear arithmetic function, e.g. the square function. We do not pursue this level of generality.
\end{rem}

\subsection*{Acknowledgments and AI disclosure}

The paper was produced with the help of various AI tools. First, the author
produced a blackboard draft containing a sketch of the main steps of the
proofs. DeepSeek-v4 was then used to translate the handwritten draft into
a \LaTeX\ file. Kimi-K3.1-pro was then used to complete the article.
GPT-6-Astra was then used to find larger values of $\delta>0$, whose proofs are too complicated to be presented in this short article.
GPT-6-Astra was also used for the revision of the article and helped to
fix typos and proof gaps. The author thank an OpenAI account holder, who is not willing to reveal itself, who generously ran the author's requests. 

We (or rather Gemini via Google search) found Idowu's paper \cite{Ido2026} during the
preparation of this article. The same exponential and logarithmic
coordinates appear there; see Theorem~4.3 and Section~7 of that
paper. 

\section{A brief outline of the proof ideas}\label{sec: idea}

We first present key ideas and leave details to a later section.

The first observation is that $2$ is a primitive root modulo $3^{K+1}$
for every $K\geq0$. Thus $2^n$, with $0\leq n<2\cdot3^K$, visits all
units of $\bZ/3^{K+1}\bZ$ exactly once. In terms of ternary expansion,
these are precisely the residues whose last digit is nonzero.
If $n\in T$, this last digit must be $1$, so $n$ is even. We therefore
write $n=2k$ and consider $4^k$.

The last digit of $4^k$ is always $1$. This is the only restriction
on its residues over a full period: the map
\[
 E_K:\bZ/3^K\bZ\longrightarrow\bZ/3^K\bZ,
 \qquad E_K(k)=\frac{4^k-1}{3}\pmod{3^K},
\]
is a permutation. We use $K$ for the number of digits remaining after
removing the final digit $1$. The inverse map
\[
 l_K:\bZ/3^K\bZ\longrightarrow\bZ/3^K\bZ
\]
can be considered as a discrete logarithm. We prove the required
arithmetic facts below.

For $K\geq0$, put
\[
 C_K=\left\{\sum_{j=0}^{K-1}\epsilon_j3^j:
                  \epsilon_j\in\{0,1\}\right\},
 \qquad C_0=\{0\}.
\]
For $K\geq1$, let $T_K=l_K(C_K)$, viewed as a subset of
$\bZ/3^K\bZ$. Also write
\[
 S=\{k\in\bN:4^k\in C\},\qquad T=\{2k:k\in S\}.
\]
If $k\in S$, removing the final digit of $4^k$ shows that
$(4^k-1)/3$ has only digits $0,1$. Hence
\begin{equation}\label{eq:suffix-inclusion}
 k\bmod3^K\in T_K\qquad\text{for every }K\geq1.
\end{equation}
In particular, $\#T_K=\#C_K=2^K$. Taking $3^K$ comparable to $N$
already gives the classical exponent $\alpha$, i.e. Narkiewicz's theorem.

Borrowing an analogy from cryptography, $T_K$ is the set of messages
$x$ whose ``encryption'' $E_K(x)$ belongs to the visibly structured
set $C_K$. The discrete logarithm should destroy this structure.
This suggests that $l_K(C_K)$ should look rather ``random''.
A useful way to make this precise is through its discrete Fourier
transform. With the normalization defined below, we prove
\begin{equation}\label{eq:sketch-decay}
 |\widehat{T_K}(\xi)|\ll3^{-\gamma K}
 \quad\text{if }3\nmid\xi,
 \qquad \gamma=\frac{2\alpha-1}{4}>0.
\end{equation}
At other nonzero frequencies divisible by powers of $3$, the same
estimate holds with $K$ replaced by a suitable smaller integer.
For $\xi\neq0$ in $\bZ/3^K\bZ$, write $s=v_3(\xi)$, so that
$3^s\mid\xi$ and $3^{s+1}\nmid\xi$, and define the conductor length
\[
 \|\xi\|_{3,K}=3^{K-s}.
\]
This is the order of the character $x\mapsto e^{2\pi i x\xi/3^K}$.
Equivalently, $\|\xi\|_{3,K}=|\xi/3^K|_3$, where $|\cdot|_3$
is the usual $3$-adic norm. Then we have
\[
 |\widehat{T_K}(\xi)|\ll \|\xi\|_{3,K}^{-\gamma},
\]
which is analogous to the Fourier decay results in Banaji--Yu \cite{BY}.

To obtain \eqref{eq:sketch-decay}, we split the $K$ digits into two
blocks of nearly equal length. Once the lower block is fixed,
$l_K$ is affine in the upper block. Cauchy--Schwarz then reduces the
problem to an $L^2$ sum of Fourier coefficients of a digit set.
We enlarge this sum to all frequencies and apply Parseval.
The saving comes from $3^\ell/4^\ell$, where $\ell=\lceil K/2\rceil$. See
Banaji--Yu \cite[Sections 3.2--3.4]{BY} for a continuous version of this
argument. In the $3$-adic setting, $l_K$ is the reduction modulo $3^K$
of $\mathcal L(x)=\log(1+3x)/\log4$ on $\bZ_3$, where both logarithms
are $3$-adic. Its derivative is
$\mathcal L'(x)=3/((1+3x)\log4)$. The logarithm series gives
\[
 \mathcal L(x+3^ry)-\mathcal L(x)-3^r\mathcal L'(x)y
 \in3^{2r+1}\bZ_3
\]
for $x,y\in\bZ_3$ and integers $r\geq0$. Thus the linearization is
exact modulo $3^K$ when $2r+1\geq K$. Lemma~\ref{lem:two-blocks}
proves the required finite congruence directly.

Finally, a standard Fourier analytic argument shows, for an interval $J$
of $L$ integers and $q=3^K\geq2L$,
\[
 \#\{k\in J:k\bmod q\in T_K\}
 \ll Lq^{\alpha-1}+q^{\alpha-\gamma}.
\]
Choosing $q\asymp L^{1/(1-\gamma)}$ balances the two terms and gives
\[
 \frac{\alpha-\gamma}{1-\gamma}
 =\frac{2\alpha+1}{5-2\alpha}=\alpha-\delta.
\]
We will make this sketch precise in the following section.

\section{The proof}

\subsection{The discrete Fourier transform}

Write $e_q(x)=\exp(2\pi i x/q)$. We use the following convention.

\begin{defn}[DFT]
Let $N\geq1$ be an integer and let $f:\bZ/3^N\bZ\to\bC$.
For $\xi\in\bZ/3^N\bZ$, define
\[
 \widehat f(\xi)=\sum_{j=0}^{3^N-1}f(j)e_{3^N}(j\xi).
\]
If $A\subset\bZ/3^N\bZ$ is nonempty, write
\[
 \widehat A(\xi)=\frac{1}{\#A}\widehat{\one_A}(\xi).
\]
\end{defn}
For sets $A$, $\hat{A}$ is not taken to be $\widehat{\one_A}$.
The normalisation makes $|\hat{A}|\leq1$. The transform depends on $N$.
We may write $\cF_N(f)$ when the
modulus is not clear from the context. Notice that the transform of
a function is unnormalized, whereas the transform of a set is
normalized to have value $1$ at zero.

For $q=3^N$, finite orthogonality gives
\begin{equation}\label{eq:orthogonality}
 \sum_{\xi=0}^{q-1}e_q(\xi x)
 =\begin{cases}q,&x\equiv0\pmod q,\\0,&x\not\equiv0\pmod q.
 \end{cases}
\end{equation}
Consequently, Fourier inversion and Parseval have the following form:
\begin{equation}\label{eq:fourier-inversion}
 f(x)=\frac1q\sum_{\xi=0}^{q-1}\widehat f(\xi)e_q(-x\xi),
 \qquad
 \sum_{\xi=0}^{q-1}|\widehat f(\xi)|^2
 =q\sum_{x=0}^{q-1}|f(x)|^2.
\end{equation}
In particular,
\begin{equation}\label{eq:set-parseval}
 \sum_{\xi=0}^{q-1}|\widehat A(\xi)|^2=\frac q{\#A}.
\end{equation}

\subsection{The discrete logarithm and the two digit blocks}

We first justify the arithmetic used in the sketch. For $m\geq0$, put
\[
 U_m=\frac{4^{3^m}-1}{3^{m+1}}.
\]
The binomial theorem gives
\[
 U_0=1,\qquad
 U_{m+1}=U_m+3^{m+1}U_m^2+3^{2m+1}U_m^3.
\]
By induction, each $U_m$ is an integer and $U_m\equiv1\pmod3$.
It follows that $4$ has order $3^K$ modulo $3^{K+1}$ for every
$K\geq1$: its order divides $3^K$ and does not divide $3^{K-1}$.
Thus its powers give all $3^K$ residues congruent to $1$ modulo $3$.
This also proves the assertion about the order of $2$, since an
exponent giving $2^n\equiv1\pmod3$ must be even. The case $K=0$
is immediate.

The map $E_K$ is therefore well defined and bijective. Its inverse
$l_K$ is characterized by
\begin{equation}\label{eq:log-definition}
 4^{l_K(y)}\equiv1+3y\pmod{3^{K+1}}.
\end{equation}
Uniqueness shows that these inverses commute with reduction:
\begin{equation}\label{eq:reduction}
 l_K(y)\equiv l_r(y\bmod3^r)\pmod{3^r}
 \qquad(1\leq r\leq K).
\end{equation}

\begin{lem}\label{lem:two-blocks}
Let $K\geq1$, $m=\lfloor K/2\rfloor$ and $\ell=\lceil K/2\rceil$.
Let $c_m$ be the inverse of $U_m$ modulo $3^\ell$.
For $u\in C_m$ and $v\in C_\ell$,
\begin{equation}\label{eq:two-blocks}
 l_K(u+3^mv)
 \equiv l_K(u)+3^mc_m(1+3u)^{-1}v\pmod{3^K}.
\end{equation}
The inverse of $1+3u$ on the right is taken modulo $3^\ell$.
\end{lem}

\begin{proof}
Since $2m+2\geq K+1$, the binomial theorem gives, for $z\geq0$,
\[
 4^{3^mz}=(1+3^{m+1}U_m)^z
 \equiv1+3^{m+1}U_mz\pmod{3^{K+1}}.
\]
Choose $z\equiv c_m(1+3u)^{-1}v\pmod{3^\ell}$. Multiplying by
$4^{l_K(u)}$ and using \eqref{eq:log-definition}, we obtain
\[
 4^{l_K(u)+3^mz}
 \equiv(1+3u)(1+3^{m+1}U_mz)
 \equiv1+3u+3^{m+1}v\pmod{3^{K+1}}.
\]
The last expression is $1+3(u+3^mv)$. The result follows from
the uniqueness in \eqref{eq:log-definition}.
\end{proof}

\subsection{The \texorpdfstring{$L^2$}{L2} bound}

We now give the main estimate. All transforms of $C_\ell$ below are
taken modulo $3^\ell$.

\begin{prop}\label{prop:l2}
Let $K,m,\ell,c_m$ be as in Lemma~\ref{lem:two-blocks}. For
$\xi\in\bZ/3^K\bZ$ with $3\nmid\xi$,
\begin{equation}\label{eq:l2-bound}
 |\widehat{T_K}(\xi)|^2
 \leq\frac1{2^m}\sum_{u\in C_m}
 \left|\widehat{C_\ell}\bigl(\xi c_m(1+3u)^{-1}\bigr)\right|^2
 \leq2\left(\frac34\right)^\ell.
\end{equation}
\end{prop}

\begin{proof}
Every element of $C_K$ is uniquely of the form $u+3^mv$, with
$u\in C_m$ and $v\in C_\ell$. Since $l_K$ is bijective,
Lemma~\ref{lem:two-blocks} gives
\begin{align*}
 \widehat{T_K}(\xi)
 &=\frac1{2^K}\sum_{u\in C_m}\sum_{v\in C_\ell}
       e_{3^K}\bigl(\xi l_K(u+3^mv)\bigr)\\
 &=\frac1{2^m}\sum_{u\in C_m}e_{3^K}\bigl(\xi l_K(u)\bigr)
       \widehat{C_\ell}\bigl(\xi c_m(1+3u)^{-1}\bigr).
\end{align*}
The first inequality in \eqref{eq:l2-bound} is Cauchy--Schwarz.

We only need to count how often a frequency on the right can occur.
Since $\xi c_m$ is a unit modulo $3^\ell$,
\[
 \xi c_m(1+3u)^{-1}\equiv\xi c_m(1+3u')^{-1}\pmod{3^\ell}
 \quad\Longleftrightarrow\quad
 u\equiv u'\pmod{3^{\ell-1}}.
\]
We have $m\geq\ell-1$. Each residue in $C_{\ell-1}$ has exactly
$2^{m-\ell+1}$ lifts to $C_m$. Thus each frequency which occurs
has weight $2^{1-\ell}$ in the average. Enlarging the frequency
sum to all residues modulo $3^\ell$ and using
\eqref{eq:set-parseval}, we obtain
\[
 |\widehat{T_K}(\xi)|^2
 \leq2^{1-\ell}\sum_{b\bmod3^\ell}
                |\widehat{C_\ell}(b)|^2
 =2^{1-\ell}\frac{3^\ell}{2^\ell}
 =2\left(\frac34\right)^\ell.
\]
This also includes $K=1$, for which $m=0$ and $C_m=\{0\}$.
\end{proof}

\begin{cor}\label{cor:decay}
Put $\gamma=(2\alpha-1)/4$. If $3\nmid\xi$, then
\[
 |\widehat{T_K}(\xi)|\leq\sqrt2\,3^{-\gamma K}.
\]
More generally, write a nonzero frequency as
$\xi=3^{K-r}b$, where $1\leq r\leq K$ and $3\nmid b$. Then
\begin{equation}\label{eq:conductor-decay}
 \widehat{T_K}(\xi)=\widehat{T_r}(b),\qquad
 |\widehat{T_K}(\xi)|\leq\sqrt2\,3^{-\gamma r}.
\end{equation}
\end{cor}

\begin{proof}
Taking square roots in \eqref{eq:l2-bound} and using $\ell\geq K/2$,
we have
\[
 |\widehat{T_K}(\xi)|
 \leq\sqrt2\left(\frac34\right)^{\ell/2}
 \leq\sqrt2\left(\frac34\right)^{K/4}
 =\sqrt2\,3^{-\gamma K}.
\]
For the second assertion, use \eqref{eq:reduction}. Each element
of $C_r$ has $2^{K-r}$ lifts to $C_K$, and consequently
\[
 \widehat{T_K}(3^{K-r}b)
 =\frac1{2^K}\sum_{y\in C_K}e_{3^r}\bigl(b l_K(y)\bigr)
 =\frac1{2^r}\sum_{y\in C_r}e_{3^r}\bigl(b l_r(y)\bigr)
 =\widehat{T_r}(b).
\]
Apply the first assertion at level $r$.
\end{proof}

\subsection{Counting in an interval}

We finish with an interval estimate. Intervals in the following
lemma consist of consecutive integers, and their length is their
number of elements.

\begin{lem}\label{lem:interval}
Let $J$ be an interval of $L\geq1$ integers and let $q=3^K\geq2L$.
Then
\begin{equation}\label{eq:interval-bound}
 \#\{k\in J:k\bmod q\in T_K\}
 \ll Lq^{\alpha-1}+q^{\alpha-\gamma},
\end{equation}
uniformly in the location of $J$.
\end{lem}

\begin{proof}
Choose an integer $z$ so that $|k-z|\leq L/2$ for all $k\in J$.
On $\bZ/q\bZ$, let $P=\{0,\dots,L-1\}$ and define\footnote{It is perhaps more intuitive to use $L^{-1}\one_{\{0,\dots,L-1\}}$ rather than the triangular weights $g$. This was the author's original choice. GPT-6-Astra suggested the usage of $g$, which has the advantage of removing a $\log q$ factor in the counting estimate. Such a factor has no effect in the final result for $\delta$.}
\[
 g(x)=\frac2L\sum_{j\in P}\one_P(x-z+j).
\]
This is a nonnegative triangular weight. If $x$ is the residue of
$k\in J$, then at least $L-|k-z|\geq L/2$ of the summands equal
$1$. Hence $g(x)\geq1$. Direct summation gives
\begin{equation}\label{eq:weight}
 \widehat g(0)=\sum_{x\bmod q}g(x)=2L,
 \qquad
 |\widehat g(\xi)|=\frac2L|D_L(\xi)|^2,
\end{equation}
where
\[
 D_L(\xi)=\sum_{j=0}^{L-1}e_q(j\xi).
\]

Let $Q$ divide $q$ and let $s$ be the remainder of $L$ on division
by $Q$. For $1\leq b<Q$, complete periods in the sum defining
$D_L((q/Q)b)$ vanish. Applying \eqref{eq:orthogonality} to the
remaining $s$ terms, and removing the zero frequency, gives
\[
 \sum_{b=1}^{Q-1}|D_L((q/Q)b)|^2=Qs-s^2\leq QL.
\]
It follows from \eqref{eq:weight} that
\begin{equation}\label{eq:weight-conductor}
 \sum_{b=1}^{Q-1}|\widehat g((q/Q)b)|\leq2Q.
\end{equation}

Since $L<q$, different integers of $J$ have different residues
modulo $q$. By positivity and Fourier inversion, the left side of
\eqref{eq:interval-bound} is at most
\begin{align*}
 \sum_{x\in T_K}g(x)
 &=\frac{2^K}{q}\sum_{\xi\bmod q}
                   \widehat g(\xi)\widehat{T_K}(-\xi)\\
 &\leq\frac{2^K}{q}\left(2L+
       \sum_{r=1}^K\sqrt2\,3^{-\gamma r}
       \sum_{\substack{1\leq b<3^r\\3\nmid b}}
                 |\widehat g(3^{K-r}b)|\right).
\end{align*}
Here we grouped the nonzero frequencies by their conductors and
used Corollary~\ref{cor:decay}. Since $2^K=q^\alpha$ and $\gamma<1$,
\eqref{eq:weight-conductor} and a geometric sum give
\[
 \sum_{x\in T_K}g(x)
 \leq\frac{2^K}{q}\left(2L+
                   2\sqrt2\sum_{r=1}^K3^{(1-\gamma)r}\right)
 \ll Lq^{\alpha-1}+q^{\alpha-\gamma}.
 \qedhere
\]
\end{proof}

\begin{proof}[Proof of Theorem~\ref{thm:main}]
Let $J$ be an interval of $L\geq1$ positive integers. Choose $q$ to
be the least power of $3$ satisfying
\[
 q\geq\max\{2L,L^{1/(1-\gamma)}\}.
\]
Then
\[
 L^{1/(1-\gamma)}\leq q\leq6L^{1/(1-\gamma)}.
\]
By \eqref{eq:suffix-inclusion} and Lemma~\ref{lem:interval},
\[
 \#(S\cap J)
 \ll Lq^{\alpha-1}+q^{\alpha-\gamma}
 \ll L^{(\alpha-\gamma)/(1-\gamma)}=L^{\alpha-\delta}.
\]
In the second inequality, use the lower bound on $q$ in the first
term and the upper bound in the second, since
$\alpha-1<0<\alpha-\gamma$.

The even integers in an interval $I$ of $N$ positive integers
correspond, under $n=2k$, to an interval of at most
$\lceil N/2\rceil\leq N$ positive integers. Since $T=\{2k:k\in S\}$,
the result follows. If $I$ contains no even integer, there is
nothing to prove.
\end{proof}

\bibliographystyle{amsplain}
\begin{samepage}

\end{samepage}
\end{document}